\documentclass[reqno]{amsart}
\usepackage{amscd,amsfonts,amssymb}
\usepackage{graphicx}
\usepackage{color}

\numberwithin{equation}{section}
\newtheorem{Theorem}{Theorem}[section]
\newtheorem{Lemma}{Lemma}[section]

\theoremstyle{definition}
\newtheorem{Definition}{Definition}[section]
\theoremstyle{remark}

\newtheorem{Proposition}{Proposition}[section]

\renewcommand{\r}{\rho}
\renewcommand{\t}{\theta}

\def\va{\varphi}
\renewcommand{\u}{{\bf u}}
\renewcommand{\H}{{\bf H}}
\newcommand{\R}{{\mathbb R}}
\newcommand{\Dv}{{\rm div}}

\newcommand{\na}{\nabla}

\newcommand{\E}{{\mathcal E}}
\newcommand{\K}{{\mathcal K}}
\newcommand{\T}{{\mathcal T}}
\def\f{\frac}
\renewcommand{\O}{\Omega}
\def\ov{\overline}

\def\hf1{^\f{1}{1-\xi^2}}

\def\be{\begin{equation}}
\def\en{\end{equation}}
\def\bs{\begin{split}}
\def\es{\end{split}}

\author{Xianpeng HU and Dehua Wang}
\address{Department of Mathematics, University of Pittsburgh,
                           Pittsburgh, PA 15260, USA.}
\email{xih15@pitt.edu}
\address{Department of Mathematics, University of Pittsburgh,
                           Pittsburgh, PA 15260, USA.}
\email{dwang@math.pitt.edu}

\title[Compactness of Solutions to Magnetohydrodynamics]
{Compactness of weak solutions to the three-dimensional compressible
magnetohydrodynamic equations}
\keywords{Three-dimensional full compressible MHD equations,
compactness, weak solutions, density-dependent viscosities.}
\subjclass{ 35Q36, 35D05, 76W05.}

\date{\today}

\begin{document}

\begin{abstract}
The compactness of  weak solutions to the
magnetohydrodynamic equations for the viscous, compressible, heat conducting
fluids is considered in both  the three-dimensional space $\R^3$ and
the three-dimensional periodic domains. The viscosities, the
heat conductivity as well as the magnetic coefficient  are allowed to depend on the density, and  may vanish on the vacuum. This paper provides a different idea from
\cite{hw2} to show the compactness of  solutions of
viscous, compressible, heat conducting magnetohydrodynamic flows, derives a
new  entropy identity, and shows that
the limit of a sequence of weak solutions is still a weak solution
to the compressible magnetohydrodynamic equations.
\end{abstract}

\maketitle

\section{Introduction}

Magnetohydrodynamics (MHD) is the theory of the macroscopic
interaction of electrically conducting fluids with  magnetic fields.
It has  a very broad range of applications.
It is of importance in connection with many engineering problems,
such as sustained plasma confinement for controlled thermonuclear
fusion, liquid-metal cooling of nuclear reactors, and
electromagnetic casting of metals. It also finds applications in
geophysics and astronomy, where one prominent example is the
so-called dynamo problem, that is, the question of the origin of the
Earth's magnetic field in its liquid metal core.

Due to their practical relevance, MHD  problems have long been
the subject of intense cross-disciplinary research, but except for
relatively simplified special cases, the rigorous mathematical and
numerical analysis of such problems remains open. In the viscous
incompressible case, MHD flow is governed by the Navier-Stokes
equations and the Maxwell equations of the magnetic field. 
For the mathematical analysis in incompressible MHD equations, see
\cite{DL, gb, MR} and the references therein. In the compressible
case, the mathematical analysis is much more complicated, due to the
oscillation of the density, the concentration of the temperature,
and the coupling interaction of hydrodynamics with the magnetic
field. The full system of the three-dimensional magnetohydrodynamic equations
in the Eulerian coordinates can be read as follows(\cite{KL, LL}):
\begin{subequations} \label{1}
\begin{align}
&\r_t +\Dv(\r\u)=0,\label{11}\\
&(\r\u)_t+\Dv\left(\r\u\otimes\u\right)+\nabla p=(\na \times
\H)\times \H+\Dv\Psi,\label{12}\\
&\E_t+\Dv(\u(\E'+p))=\Dv((\u\times\H)\times\H+\nu\H\times(\nabla\times\H)+\u\Psi+\kappa\nabla\theta),\label{13}\\
&\H_t-\nabla\times(\u\times\H)=-\nabla\times(\nu\nabla\times\H),\qquad
\Dv\H=0,\label{14}
\end{align}
\end{subequations}
where $\Psi=2\mu D(\u)+\lambda\,\Dv\u\,\mathbf{I}$ with
$3\lambda+2\mu\ge 0$ and $D(\u)=\f12(\nabla\u+(\nabla\u)^\top)$
denotes the strain rate tensor; $\r$ denotes the density, $\u\in
\R^3$ the velocity, $\H\in \R^3$ the magnetic field, and $\theta$
the temperature; $\E$ is the total energy given by
$$\E=\r\left(e+\f{1}{2}|\u|^2\right)+ \f{1}{2}|\H|^2 \textrm{ and }
\E'=\r\left(e+\f{1}{2}|\u|^2\right),$$ with $e$ the internal energy,
$\f{1}{2}\r|\u|^2$ the kinetic energy, and $\f{1}{2}|\H|^2$ the
magnetic energy. The equations of state $p=p(\r,\theta)$,
$e=e(\r,\theta)$ relate the pressure $p$ and the internal energy $e$
to the density and the temperature of the flow; $\mathbf{I}$ is the
$3\times 3$ identity matrix, and $(\nabla\u)^\top$ is the transpose
of the matrix $\nabla\u$. $\nu(\r, \t)$ is the magnetic field
coefficient, $\kappa=\kappa(\r, \theta)$ is the heat conductivity.
In general, equations \eqref{11}, \eqref{12}, \eqref{13} denote the
conservations of mass,
 momentum, and energy, respectively.
The equation \eqref{14} is called the induction equation, and the
electric field can be written in terms of the magnetic field $\H$
and the velocity $\u$,
\begin{equation*} 
{\bf E}=\nu\nabla\times\H - \u\times\H.
\end{equation*}
Although the electric field ${\bf E}$ does not appear in the MHD
system \eqref{11}-\eqref{14}, it is indeed induced according to
the above relation by the moving conductive flow in the magnetic
field.

In this paper, we are interested in the compactness of weak
solutions to the compressible MHD equations \eqref{1} both in the
three-dimensional space $\R^3$ and in the three-dimensional periodic
domains. As it is well-known, the motivation of considering the
compactness of weak solutions is to show the existence of weak
solutions and the stability of weak solutions of nonlinear problems.
In the literature, there have been a lot of studies on MHD by
physicists and mathematicians because of its physical importance,
complexity, rich phenomena, and mathematical challenges; see
\cite{gw, gw2, f6, FJN, h1, HT, hw1, hw2, sm, LL, w1} and the
references cited therein. For instance, the smooth global solution
near the constant state in one-dimensional case is investigated in
\cite{sm}. However, many fundamental problems for compressible MHD
with large, discontinuous initial data are still open.

Positive results on the existence of weak
solutions with large, discontinuous data for compressible MHD equations have been obtained
recently in \cite{hw1, hw2}, specially in \cite{hw2} for full
compressible MHD equations with temperature-dependent viscosities.
More precisely, it was shown in \cite{hw2}, under certain
structural hypotheses imposed on the pressure $p$ and the heat
conductivity coefficient $\kappa$, that the full compressible MHD
system admits at least a global-in-time variational solution for
large initial data. Those solutions satisfy the equations
\eqref{11}, \eqref{12}, \eqref{14} in the sense of distributions
while the thermal energy equation \eqref{13} is being replaced by
two inequalities to be accordance with the second law of thermodynamics.
 This approach is in the spirit of the concept of \textit{
weak solutions with a defect measure } introduced by several
authors in different contexts, see \cite{dl}.
However, in order to obtain the estimates on the gradient of the
velocity, the works in \cite{hw1, hw2} do rely strongly on the
assumption that the shear viscosity $\mu$ is bounded below by a
positive constant.

Our aim, in this paper, is to show the compactness of weak solutions
to the full compressible MHD equations with viscosity coefficients
vanishing on the vacuum both in the three-dimensional space $\R^3$
and in the three-dimensional periodic domains. Although the periodic
case does not correspond to a physical configuration, its
mathematical treatment is technically easier, while it retains the
main mathematical difficulties of the problem of the flow. More
importantly, in our context, the viscosities $\mu, \lambda$, the heat
conductivity $\kappa$, and the magnetic coefficient $\nu$ can be
allowed to depend on the density $\r$ and the temperature $\t$ of
the flow. We remark that the similar problems for the compressible
Navier-Stokes equations have been studied in \cite{bd1, bd2, bdg,
mv}. Comparing with those works on the compressible Navier-Stokes
equations, we will encounter extra difficulties in studying the
compressible MHD equations. More precisely, besides the possible
oscillation of the density and the concentration of the temperature,
the appearance of the magnetic field and the coupling effect between
the hydrodynamic flow and the magnetic field should also been taken
into consideration. We remark that a simultaneous independent work similar to this paper for the case of the two-dimensional periodic domain was done in \cite{Sart}.

The novelty of this paper is to provide a new method to deal with
the vanishing viscosities for compressible MHD flows. The loss of
positivity of the viscosity coefficients implies that there is no
hope to obtain directly the uniform bound on the gradient of the
velocity. It is well known that the main difficulty, in proving
the compactness of weak solutions of compressible MHD equations,
is to pass to the limit for the nonlinear terms. In our context,
the new kind of entropy equality will provide the estimates on the
gradient of the density, which  makes the nonlinear terms much
easier to be dealt with and also give rise to a new estimate of
$\r\u^2$ in a functional space better than $L^\infty([0,T];
L^1(\O))$. In other words, although the regularity on the velocity that we can
get directly is much lower, the regularity on the density in our
context is much higher. To achieve this aim, the entropy equation
\eqref{29} and the thermal equation \eqref{41} need to be taken
into consideration. But, unfortunately, the case with constant
viscosity coefficients is excluded from our setting and the
extension to the general case in which the viscosity coefficients
depend on both the density and the temperature seems also out of
the reach of our present work.

We organize the rest of this paper as follows. In Section 2, we
will give the hypotheses in detail, introduce the definition of
weak solutions, and state our compactness result (Theorem
\ref{mt}). In Section 3, we will derive the \textit{ a priori}
estimates and a new kind of the entropy identity. In Section 4,
some auxiliary integrability lemmas are showed. Finally, we will
finish the proof of Theorem \ref{mt} in Section 5 using
Aubin-Lions Lemma.

\bigskip

\section{ Assumptions and the Main Result}
To our best knowledge, the rigorous mathematical analysis for
compressible flows is beyond the available mathematical framework.
Hence, we need add some restrictions to viscosity coefficients
$\mu$, $\lambda$, the heat conductivity $\kappa$ and the magnetic
coefficient $\nu$.
\subsection{Assumptions}
To begin with, we assume that $\mu(\r)$ and $\lambda(\r)$ are two
$C^1(0,\infty)$ functions satisfying
\begin{equation}\label{31}
\lambda(\r)=2(\r\mu'(\r)-\mu(\r)).
\end{equation}
As seen later on, this relation is fundamental to get higher
regularity on the density. More precisely, with the help of this
relation, we can show a new kind of entropy equality, which then
gives the uniform bound on the gradient of the density. Next, due
to our technical restrictions, we will need the following
constraints: $\mu(0)=0$ and there exist positive constants $c_0$,
$c_1$, $A$ and $m>1$, $2/3<\beta<1$ such that
\begin{equation}\label{32}
\begin{cases}
&\textrm{ for all } s<A,\quad c_0 s^{\beta-1}\le\mu'(s)\le
\f{s^{\beta-1}}{c_0} \textrm{ and }3\lambda(s)+2\mu(s)\ge c_0 s^\beta,\\
&\textrm{ for all } s\ge A,\quad c_1 s^{m-1}\le\mu'(s)\le
\f{s^{m-1}}{c_1} \textrm{ and } c_1 s^m\le3\lambda(s)+2\mu(s)\le
\f{s^m}{c_1}.
\end{cases}
\end{equation}
Observing that the assumption \eqref{32} implies that $\mu'(\r)>0$
for $\r>0$.

The heat conductivity coefficient $\kappa$ is assumed to satisfy
\begin{equation}\label{33}
\kappa(\r, \t)=\kappa_0(\r, \t)(\r+1)(\t^a+1),
\end{equation}
where $a\ge 2$, and $\kappa_0$ is a $C^0(R_{+}^2)$ function
satisfying for all positive $\r$ and $\t$,
$$c_2\le\kappa_0(\r, \t)\le\f{1}{c_2},$$
for some positive constant $c_2$.

For the pressure, we assume that the equations of state are of ideal
polytropic gas type
\begin{equation}\label{34}
p=\r\t+p_e(\r), \quad e=c_\nu\t+P_e(\r),
\end{equation}
with $P_e(\r)=\int_1^\r p_e(\xi)/\xi^2\mathrm{d}\xi$. We also
require that $p_e(\r)$ satisfies
\begin{equation}\label{35}
\begin{cases}
&c_3\r^{-l-1}\le p_e'(\r) \le\f{1}{c_3}\r^{-l-1}, \textrm{ if } 0\le\r<A_0,\\
&p_e'(\r)\le c_4\r^{k-1}, \textrm{ if } \r>A_0,
\end{cases}
\end{equation}
for some $A_0>0$, $c_3>0$, $c_4>0$, $l>\f{2\beta(3m-2)}{m-1}-1$,
and $k\le \left(m-\f{1}{2}\right)\f{5(l+1)-6\beta}{l+1-\beta}$.

For the magnetic coefficient, we need the following assumption:
\begin{equation}\label{36}
\nu(\r, \t)\ge c_5\f{\t}{\r}\textrm{ on }\{\r>0\},\quad
\f{1}{c_6}\ge\nu(\r, \t)\ge c_6,
\end{equation}
for some $c_5> 0$, $c_6>0$.

\subsection{Main Result}
Before we state the compactness result, we need to specify the
definition of weak solutions which we will address. It is
necessary to require that the weak solutions should satisfy the
natural energy estimates and from the viewpoint of physics, the
conservation laws on  mass,  momentum and energy also
should be satisfied at least in the sense of distributions. Based
on those considerations, the definition of reasonable
global-in-time weak solutions goes as follows.
\begin{Definition} A vector $(\r, \u, \t, \H)$ is said to be a
global-in-time weak solution to the full compressible MHD system
\eqref{11}-\eqref{14}, if and only if for any positive number $T$,
the following conditions are satisfied:
\begin{itemize}
\item \begin{equation*}
\begin{split}
\r\E\in L^\infty([0,T]; L^1(\O)),\quad \r|\u|^2\in L^\infty([0,T];
L^1(\O)),\\ \f{\nabla\mu(\r)}{\sqrt{\r}}\in L^\infty([0,T];
L^2(\O)),\quad
(\r^{\beta/2}+\r^{m/2})\nabla\u\in L^2([0,T]; L^2(\O)),\\
(1+\sqrt{\r})\nabla\t^{a/2}\in L^2([0,T]; L^2(\O)),\quad \quad
(1+\sqrt{\r})\f{\nabla\t}{\t}\in L^2([0,T]; L^2(\O)),
\end{split}
\end{equation*}
for $a\ge 2$. Moreover, for large enough $s>0$, we have
$$\r,\quad \r\u, \quad \r\E,\quad \H \in C([0,T]; H^{-s}(\O)),$$ and
hence, the initial data are satisfied in the sense of
distributions (denoted  by $\mathcal{D}'(\O\times(0,T))$.

\item The equations \eqref{11}-\eqref{14} are satisfied in the sense
of distributions.
\end{itemize}
\end{Definition}

Now our compactness result can be read as follows:
\begin{Theorem}\label{mt}
Let $\O$ be either the three-dimensional periodic domains or the
three-dimensional space $\R^3$. Assume that $\mu$, $\lambda$, $\nu$,
$\kappa$ are $C^1[0,\infty)$ functions satisfying the hypotheses
\eqref{31}-\eqref{36}. Let $\{(\r_n, \u_n, \t_n,
\H_n)\}_{n=1}^{\infty}$ be a sequence of weak solutions of
\eqref{11}-\eqref{14} satisfying the entropy equation \eqref{29} and
the thermal energy equation \eqref{41} with the initial data
\begin{equation*}
\begin{split}
\r_n(x,0)=\r_{0,n}(x),\quad \u_n(x,0)=\u_{0,n}(x),\\
\t_n(x,0)=\t_{0,n}(x),\quad \H_n(x,0)=\H_{0,n}(x),
\end{split}
\end{equation*}
where $\r_{0,n}$, $\u_{0,n}$, $\t_{0,n}$, $\H_{0,n}$ satisfy
\begin{equation}\label{initial}
\begin{cases}
&\r_{0,n}\ge 0,\quad \r_{0,n}\rightarrow \r_0 \textrm{ in } L^1(\O),\\
&\r_{0,n}\ln\r_{0,n}\rightarrow\r_0\ln\r_0\in L^1(\O), \quad \r_{0,n}\ln\t_{0,n}\rightarrow\r_0\ln\t_0\in L^1(\O),\\
&\r_{0,n}|\u_{0,n}|^2\rightarrow \r_0|\u_0|^2\textrm{ in } L^1(\O),\\
&\t_{0,n}>0,\quad \H_{0,n}\rightarrow \H_0 \textrm{ in } L^2(\O),\\
&\r_{0,n} e_{0,n}\rightarrow\r_0 e_0\in L^1(\O),\quad
\f{\nabla\mu(\r_0)}{\sqrt{\r_0}}\rightarrow\f{\nabla\mu(\r_0)}{\sqrt{\r_0}}\in
L^2(\O).
\end{cases}
\end{equation}
Assume also that the sequence of the magnetic field
$\{\H_n\}_{n=1}^{\infty}$ is uniformly bounded in
$L^{\infty}(\O\times(0,T))$. Then, up to a subsequence, $\{(\r_n,
\u_n, \t_n, \H_n)\}_{n=1}^{\infty}$ converges to a global-in-time
weak solution $(\r, \u, \t, \H)$ with initial data $(\r_0, \u_0,
\t_0, \H_0)$. More precisely,
\begin{equation*}
\begin{cases}
&\{\r_n\}_{n=1}^{\infty} \textrm{ converges strongly to } \r
\textrm{ in } C([0,T];
L^p(\O')) \textrm{ for all } 1\le p<6m-3;\\
&\{\u_n\}_{n=1}^{\infty} \textrm{ converges weakly to } \u
\textrm{ in }
L^{q_1}([0,T]; W^{1,q_3}(\O')) \textrm{ for } q_1>5/3, q_3>15/8;\\
&\{\t_n\}_{n=1}^{\infty} \textrm{ converges strongly to } \t
\textrm{ in } L^p([0,T]; L^q(\O')) \textrm{ for all } p<a \textrm{
and }
q<3a \textrm{ with }\\
&\qquad\f{1}{q}=\f{a-p}{p(a-1)r}+\f{p-1}{3p(a-1)}, \textrm{ for all } r<3/2;\\
&\{\H_n\}_{n=1}^{\infty} \textrm{ converges weakly to } \H \textrm{
in } L^2([0,T]; H^1(\O))\cap C([0,T]; L^2_{weak}(\O)),
\end{cases}
\end{equation*}
where $\O'$ is any sufficiently smooth and compact subset of $\O$.
\end{Theorem}


\bigskip

\section{ Energy Estimates and the Entropy Inequality}
In this section, we dedicate to the well-known \textit{a priori}
estimates and a new kind of the entropy equality on weak solutions
of the compressible MHD system \eqref{11}-\eqref{14}. To begin
with, from the total energy equation \eqref{13}, the physical
energy inequality holds
\begin{equation}\label{21}
\begin{split}
E(t):=\int_{\O}\r&\left(e+\f{1}{2}|\u|^2\right)(t,x)+\f{1}{2}|\H|^2(t,x)\,
dx\\&
\le\int_{\O}\r_0\left(e_0+\f{1}{2}|\u_0|^2\right)+\f{1}{2}|\H_0|^2
\, dx:=E(0).
\end{split}
\end{equation}

As shown in \cite{hw1, hw2}, the  energy estimate \eqref{21} alone
is not sufficient to build up a reasonable compactness theory of
weak solutions to compressible MHD equations in the sense of
distributions since we can not obtain any \textit{a priori}
estimate on the dissipation about the viscous stress and the
gradient of the magnetic field. Comparing with \textit{a priori}
estimates for isentropic
 cases (see \cite{hw1}), this is a major difference, because in the isentropic case, the viscous dissipation
 naturally provides a
$H^1$ bound in spatial variables on the velocity $\u$. To establish
the compactness theory in our \textit{new} framework, the following
calculation is crucial:
\begin{Lemma}\label{l1}
\begin{equation}\label{22}
\begin{split}
\f{1}{2}\f{d}{dt}\int_\O&\left(\r|\u|^2+|\H|^2\right)\, dx+\int_\O
2\mu(\r)D(\u):D(\u)\, dx
\\&+\int_\O \lambda(\r)|\Dv\u|^2\, dx+\int_\O\nu|\nabla\times\H|^2\, dx=\int_\O p(\r, \t)\Dv\u\, dx,
\end{split}
\end{equation}
and
\begin{equation}\label{23}
\begin{split}
\f{1}{2}&\f{d}{dt}\int_\O \left(\r|\u+2\nabla\va(\r)|^2+|\H|^2\right)\, dx+\int_\O 2\mu(\r)A(\u):A(\u)\, dx
+\int_\O\nu|\nabla\times\H|^2\, dx\\
&=\int_\O p(\r, \t)\Dv\u\, dx-2\int_\O \nabla p(\r,
\t)\cdot\nabla\va(\r)\, dx+2\int_\O
(\nabla\times\H)\times\H\cdot\f{ \nabla\mu(\r)}{\r} \, dx,
\end{split}
\end{equation}
where $A(\u)=(\nabla\u-\nabla\u^\top)/2$ denotes the skew symmetric
part of $\nabla\u$, and $\va'(x)=\mu'(x)/x$ for all $x>0$. The
notation $A:B$ denotes the dot product between two $n\times n$
matrices $A$ and $B$.
\end{Lemma}
\begin{proof}
The energy equality \eqref{22} is classical, and can be shown by multiplying the momentum equation \eqref{12}
 by $\u$, the mass conservation equation \eqref{11} by $|\u|^2/2$, and the magnetic equation \eqref{14} by $\H$,
  then summing them together. Here we used the following identity:
\begin{equation*}
\int_{\O}(\na \times \H)\times \H\cdot\u \,
dx=-\int_{\O}(\nabla\times(\u\times\H))\cdot\H \, dx.
\end{equation*}

Now, we turn to show the equality \eqref{23}. The idea is taken
from \cite{bd1}, and the argument goes as follows. From the mass
conservation equation, we deduce that
$$\partial_t\va(\r)+\u\cdot\nabla\va(\r)+\va'(\r)\r\Dv\u=0.$$
This gives, differentiating this equation with respect to the space variable $x_i$, $i=1,2,3$, noting $x=(x_1,x_2,x_3)$,
$$\partial_t\partial_i\va(\r)+(\u\cdot\nabla)\partial_i\va(\r)+(\partial_i\u\cdot\nabla)\va(\r)+\partial_i(\va'
(\r)\r\Dv\u)=0.$$ Let us multiply this equation by
$\r\partial_i\va(\r)$ and sum over $i$, by using the mass
equation, then one can deduce
\begin{equation}\label{24}
\f{1}{2}\f{d}{dt}\int_\O \r|\nabla\va(\r)|^2\,
dx+\int_\O\r\nabla\va(\r)\otimes\nabla\va(\r):\nabla\u \, dx+
\int_\O \nabla(\va'(\r)\r\Dv\u)\cdot\nabla\mu(\r)\, dx=0.
\end{equation}
Multiplying the momentum equation by $\nabla\mu(\r)/\r$, we get
\begin{equation}\label{25}
\begin{split}
&\int_\O(\partial_t\u+\u\cdot\nabla\u)\cdot\nabla\mu(\r)\,
dx+2\int_\O \mu(\r)D(\u):\left(\f{\nabla
\nabla\mu(\r)}{\r}-\f{\nabla\mu(\r)\otimes
\nabla\r}{\r^2}\right)\, dx\\
&\quad +\int_\O\nabla p(\r,
\t)\cdot\f{\nabla\r}{\r}\mu'(\r)\, dx
+2\int_\O\nabla((\mu(\r)-\mu'(\r)\r)\Dv\u)\cdot\f{\nabla\mu(\r)}{\r}\,
dx\\&=\int_\O (\nabla\times \H)\times\H\cdot\f{\nabla\mu(\r)}{\r}
\, dx.
\end{split}
\end{equation}

Integrating by parts, equation \eqref{24} can be rewritten under the form
\begin{equation}\label{26}
\begin{split}
&\f{1}{2}\f{d}{dt}\int_\O\r|\nabla\va(\r)|^2\, dx+\int_\O
\f{\mu(\r)\nabla\mu(\r)\cdot\nabla\u
\cdot\nabla\r}{\r^2}\mathrm{d}
-\int_\O\f{\mu(\r)\nabla\u:\nabla\nabla\mu(\r)}{\r}\,
dx\\&-\int_\O\f{\mu(\r)\nabla\Dv\u\cdot \nabla\mu(\r)}{\r}\, dx
+\int_\O \nabla(\va'(\r)\r\Dv\u)\cdot\nabla\mu(\r)\, dx=0.
\end{split}
\end{equation}

Adding equation \eqref{25} to equation \eqref{26} multiplied by 2, we get
\begin{equation}\label{27}
\begin{split}
&\int_\O(\partial_t\u+\u\cdot\nabla\u)\cdot\nabla\mu(\r)\,
dx-2\int_\O\f{\mu(\r)\nabla\Dv\u
\cdot\nabla\mu(\r)}{\r}\, dx\\
&\quad+\f{1}{2}\f{d}{dt}\int_\O 2\r|\nabla\va(\r)|^2\, dx+2\int_\O
\nabla((\mu(\r)-\mu'(\r)\r)\Dv
\u)\cdot\f{\nabla\mu(\r)}{\r}\, dx\\
&\quad+2\int_\O\nabla(\mu'(\r)\Dv\u)\cdot\nabla\mu(\r)\, dx+\int_\O
\nabla p(\r, \t)\cdot\f{\nabla \r}{\r}\mu'(\r)\, dx
\\&=\int_\O (\nabla\times\H)\times\H\cdot\f{\nabla\mu(\r)}{\r}
\, dx.
\end{split}
\end{equation}
By splitting the terms involving $\Dv\u$ and by summing them, we get
\begin{equation}\label{28}
\begin{split}
&\int_\O(\partial_t\u+u\cdot\nabla\u)\cdot\nabla\mu(\r)\,
dx+\f{1}{2}\f{d}{dt}\int_\O 2 \r|\nabla\va(\r)|^2\, dx+ \int_\O
\nabla p(\r, \t)\cdot\f{\nabla\r}{\r}\mu'(\r)\, dx\\&=\int_\O
(\nabla\times\H) \times\H\cdot\f{\nabla\mu(\r)}{\r} \, dx.
\end{split}
\end{equation}
But as for the first term in \eqref{28}, we can calculate
$$\int_\O(\partial_t\u+u\cdot\nabla\u)\cdot\nabla\mu(\r)\, dx=\f{d}{dt}\int_\O\u\cdot
\nabla\mu(\r)\, dx-\int_\O\u\cdot\nabla\partial _t\mu(\r)\,
dx+\int_\O(\u\cdot\nabla\u)\cdot\nabla\mu(\r)\, dx.$$
By using now the mass equation and integrating by parts the last two terms, this gives
\begin{equation*}
\begin{split}
\int_\O(\partial_t\u+u\cdot\nabla\u)\cdot\nabla\mu(\r)\,
dx&=\f{d}{dt}\int_\O\u\cdot \nabla\mu(\r)\,
dx-\int_\O\mu'(\r)\Dv(\r\u)\Dv\u \,
dx\\&\quad-\int_\O\mu(\r)\u\cdot\nabla\Dv\u\,
dx-\int_\O\mu(\r)\partial_i\u_j
\partial_j\u_i\, dx.
\end{split}
\end{equation*}
Integrating by parts in the third term, we get
\begin{equation*}
\begin{split}
\int_\O(\partial_t\u+u\cdot\nabla\u)\cdot\nabla\mu(\r)\,
dx=&\f{d}{dt}\int_\O\u\cdot \nabla\mu(\r)\, dx-\int_\O(\r\mu'(\r)-
\mu(\r))|\Dv\u|^2\, dx\\
&-\int_\O\mu(\r)\partial_i\u_j\partial_j\u_idx.
\end{split}
\end{equation*}
Adding the above identity with \eqref{28}, we get the following equality
\begin{equation*}
\begin{split}
&\f{d}{dt}\int_\O\u\cdot\nabla\mu(\r)dx-\int_\O(\r\mu'(\r)-\mu(\r))|\Dv\u|^2dx-
\int_\O\mu(\r)\partial_i\u_j\partial_j\u_idx\\&\quad+\f{1}{2}\f{d}{dt}\int_\O
2\r|\nabla\va (\r)|^2dx+\int_\O \nabla p(\r, \t)\cdot
\f{\nabla\r}{\r}\mu'(\r)dx\\&=\int_\O
(\nabla\times\H)\times\H\cdot\f{\nabla\mu(\r)}{\r} dx.
\end{split}
\end{equation*}
Adding this last equation multiplied by 2 to the energy estimate \eqref{22} gives \eqref{23}.
\end{proof}

Next, we introduce the concept of the entropy $s(\r, \t)$ which
satisfies the entropy equation (\cite{f6, hw2, KL, LL}):
\begin{equation}\label{29}
\partial_t(\r s)+\Dv(\r
s\u)+\Dv\left(\f{\kappa(\t)}{\t}\right)
=\f{1}{\t}(\nu|\nabla\times\H|^2+\Psi:\nabla\u)-\f{\kappa(\t)
|\nabla\t|^2}{\t^2},
\end{equation}
with $s(\r,\t)=c_{\upsilon}\ln\t-\ln\r$, where $c_\upsilon$ is a
positive constant denoting the specific heat at constant volume.
The entropy equation is useful in compressible flows, because it
provides naturally the estimates in the gradient of the
temperature. More precisely, integrating \eqref{29} over
$\O\times(0,t)$, the following proposition is verified:
\begin{Proposition}\label{p31}
Assume that $\r_0 s_0\in L^1(\O)$. Then, for all $t\ge 0$, one has:
\begin{equation}\label{210}
\int_0^t\!\!\!\!\int_\O\f{1}{\t}(\nu|\nabla\times\H|^2+\Psi:\nabla\u)+\f{\kappa(\r,
\t)|\nabla \t|^2}{\t^2}dx\mathrm{d}t\le \int_\O \r s+|\r_0 s_0|dx.
\end{equation}
\end{Proposition}

\medskip
Observing that
$$\r s\le c_\upsilon\r\t-\r\ln\r.$$
the first term on the right-hand side of \eqref{210} can be
estimated by
$$\int_\O \r sdx\le\int_\O
c_\upsilon\r\t dx-\int_\O\r\ln\r dx.$$
Multiplying \eqref{11} by $1+\ln\r$, we get
$$\partial_t(\r\ln\r)+\Dv(\r\u\ln\r)+\r\Dv\u=0.$$ Thus, we have
\begin{equation*}
\int_\O\r\ln\r dx=\int_\O
\r_0\ln\r_0dx+\int_0^t\!\!\!\!\int_\O\r\Dv\u \,dx\mathrm{d}s.
\end{equation*}
Therefore, the right-hand side of \eqref{210} can be estimated by
\begin{equation}\label{211}
\begin{split}
\int_\O \r s\,dx&\le\int_\O c_\upsilon\r\t\,dx+\int_\O
|\r_0\ln\r_0|\,dx+\int_0^t\!\!\!\!\int_\O\r|\Dv\u|\,dx\mathrm{d}s\\&\le\int_\O
c_\upsilon\r\t\,dx+\int_\O |\r_0\ln\r_0|\,dx+\int_0^t\!\!\!\!\int_\O
\f{\sqrt{\r}}{\sqrt{3\lambda+2\mu}}\f{\sqrt{3\lambda+2\mu}}{\sqrt{\t}}|\Dv\u|\sqrt{\r\t}\,dx\mathrm{d}t,
\end{split}
\end{equation}
and, then  using the classical Young's inequality, the bound
of $\r\t$ in $L^\infty([0,T]; L^1(\O))$, and the assumption
\eqref{32} that ensures that $ s\mapsto s/3\lambda(s)+2\mu(s)$
belongs to $L^\infty(R_{+})$, we conclude that the terms on the left-hand side of \eqref{210} is bounded in $L^1(\O\times(0,T))$.

Hence, if $\r_0 s_0$ and $\r_0\ln\r_0$ belong to $L^1(\O)$, then
the components of the following four quantities
$\sqrt{3\lambda+2\mu}|\Dv\u|/\sqrt{\t}$,
$\sqrt{\mu}D(\u)/\sqrt{\t}$, $\sqrt{\nu}\nabla\times\H/\sqrt{\t}$,
$(\sqrt{\r}+1)\nabla\t^{\f{a}{2}}$ and $(\sqrt{\r}+1)\nabla\ln\t$
are bounded in $L^2(\O\times(0,T))$. We note that the last two
bounds involving the temperature gradient provide the following
useful estimates:
\begin{equation}\label{212}
(\sqrt{\r}+1)\nabla\t^\alpha\in L^2(\O\times(0,T)), \textrm{ for
all } \alpha \textrm{ such that } 0\le\alpha\le a/2.
\end{equation}

In order to get enough \textit{a priori} estimates from Lemma
\ref{l1} and the initial condition \eqref{initial}, we have to
control the following terms:
$$\int_\O p\Dv\u\,dx, \quad \int_\O\nabla
p\cdot\nabla\va(\r)\,dx,\quad
\int_\O(\nabla\times\H)\times\H\cdot\f{\nabla\mu(\r)}{\r}\,dx.$$ To
this end, the following estimates are useful:
\begin{Lemma}\label{l2}
Let $\O$ be the three-dimensional periodic box or the whole space
$\R^3$. For all $\r$ satisfying $\r^{-1/2}\nabla\mu(\r)\in
L^2(\O)$, one has
\begin{equation*}
\begin{cases}
&\|\r^{m-1/2}\chi_{\{\r>2A\}}\|_{L^6(\O)}\le
c\|\f{\nabla\mu(\r)}{\sqrt{\r}}\|_{L^2(\O)},\\
&\|\r^{\beta-1/2}\chi_{\{\r\le A/2\}}\|_{L^6(\O)}\le
c\|\f{\nabla\mu(\r)}{\sqrt{\r}}\|_{L^2(\O)},
\end{cases}
\end{equation*}
for some positive constant $c$, where $A$ is from \eqref{32}, and
$\chi$ is the characteristic function.
\end{Lemma}
\begin{proof}
Let us consider the function $\eta=\alpha\xi^{m-1/2}$, where $\xi:
(0, \infty)\rightarrow (0, \infty)$ is a smooth increasing
function such that $\xi(s)=s$ for $s>2A$ and $\xi(s)=0$ for $s<A$
and $\alpha$ is a positive constant. By hypothesis \eqref{32}, we
can choose $c$ such that $\eta'(s)\le c\mu'(s)/\sqrt{s}$ for all
$s>0$. Using Sobolev's inequality, we have
$$\|\eta(\r)\|_{L^6(\O)}\leq c\|\nabla\eta(\r)\|_{L^2(\O)}\leq c\left\|\f{\nabla\mu(\r)}{\sqrt{\r}}\right\|_{L^2(\O)}.$$
The left-hand side of above inequality is bigger than
$\|\r^{m-1/2}\chi_{\{\r>2A\}}\|_{L^6(\O)}$. This implies
$$\|\r^{m-1/2}\chi_{\{\r>2A\}}\|_{L^6(\O)}\le
c\left\|\f{\nabla\mu(\r)}{\sqrt{\r}}\right\|_{L^2(\O)}.$$

To prove the second part, a similar approach can be applied.
Indeed, choosing the function $\eta=\alpha\xi^{\beta-1/2}$ such
that $|\eta'(s)|\le c\mu'(s)/\sqrt{s}$ for all $s>0$, where $\xi:
(0, \infty)\rightarrow (0, \infty)$ is a smooth positive function
such that $\xi(s)=0$ for $s>A$ and $\xi(s)=s$ for $s<A/2$. By
Sobolev's inequality, we have
$$\|\eta(\r)\|_{L^6(\O)}\leq c\|\nabla\eta(\r)\|_{L^2(\O)}\leq c\left\|\f{\nabla\mu(\r)}{\sqrt{\r}}\right\|_{L^2(\O)}.$$
The left-hand side of above inequality is bigger than
$\|\r^{\beta-1/2}\chi_{\{\r<A/2\}}\|_{L^6(\O)}$. This implies
$$\|\r^{\beta-1/2}\chi_{\{\r<A/2\}}\|_{L^6(\O)}\le
c\left\|\f{\nabla\mu(\r)}{\sqrt{\r}}\right\|_{L^2(\O)}.$$
\end{proof}

From Lemma \ref{l2}, we know $\r\in L^\infty([0,T];
L^{6m-3}(\O'))$ for any bounded subset $\O'$ of $\O$.

\begin{Lemma}[\bf{The control of $\int_\O p\,\Dv\u\,dx$}]\label{l3}
\begin{equation*}
\begin{split}
\int_\O p\,\Dv\u\,dx &\le -\f{d}{dt}\int_\O \r P_e(\r)\,dx+
\varepsilon\|\sqrt{3\lambda+2\mu}\Dv\u\|_{L^2}^2\\
&\quad +c_\varepsilon\left(\|\r\t\|^2_{L^1}+\|\t\|^2_{L^6}+
\|\t\|^2_{L^3} \left\|\f{\nabla\mu(\r)}{\sqrt{\r}}\right\|^2_{L^2}\right),
\end{split}
\end{equation*}
for all positive $\varepsilon$.
\end{Lemma}
\begin{proof}
By the continuity equation in the renormalized sense, we have
\begin{equation}\label{213}
\begin{split}
\int_\O p\,\Dv\u\,dx&=\int_\O p_e(\r)\Dv\u\,dx+\int_\O
\r\t\Dv\u\,dx\\&=-\f{d}{dt}\int_\O \r P_e(\r)\,dx+\int_\O
\r\t\Dv\u\,dx.
\end{split}
\end{equation}
For the second term on the right-hand side of \eqref{213}, we can
estimate
\begin{equation*}
\begin{split}
\left|\int_\O \r\t\Dv\u\,dx\right|
&\le\|\sqrt{3\lambda+2\mu}\Dv\u\|_{L^2}\\&\quad\times(\|\r\t\chi_{\{\r<A\}}/\sqrt{3\lambda+2\mu}\|_{L^2}
+\|\r\t\chi_{\{\r\ge A\}}/\sqrt{3\lambda+2\mu}\|_{L^2})\\&\le
c\|\sqrt{3\lambda+2\mu}\Dv\u\|_{L^2}\\&\quad\times(\|\r^{2/5}\t\|_{L^2}\|
\r^{(6-5\beta)/10}\chi_{\{\r<A\}}\|_{L^\infty}
+A^{-m/2}\|\r\chi_{\{\r\ge A\}}\|_{L^6}\|\t\|_{L^3})\\&\le
c\|\sqrt{3\lambda+2\mu}\Dv\u\|_{L^2}\\&\quad\times(\|\r\t\|^{2/5}_{L^1}\|
\t\|^{3/5}_{L^6}A^{(6-5\beta)/10} +A^{(3-3m)/2}\|\r\chi_{\{\r\ge
A\}}\|^{m-1/2}_{L^{6m-3}}\|\t\|_{L^3}).
\end{split}
\end{equation*}
Thus, in view of Lemma \ref{l2} and Young's inequality, we have
\begin{equation*}
\begin{split}
\left|\int_\O \r\t\Dv\u\,dx\right|&\le
c\|\sqrt{3\lambda+2\mu}\Dv\u\|_{L^2}\\&\quad\times\left(\|\r\t\|_{L^1}+\|\t\|_{L^6}+\|\t\|_{L^3}
\left\|\f{\nabla\mu(\r)}{\sqrt{\r}}\right\|_{L^2}\right)\\
&\le\varepsilon\|\sqrt{3\lambda+2\mu}\Dv\u\|_{L^2}^2\\&\quad+c_\varepsilon\left(\|\r\t\|^2_{L^1}+\|\t\|^2_{L^6}+
\|\t\|^2_{L^3} \left\|\f{\nabla\mu(\r)}{\sqrt{\r}}\right\|^2_{L^2}\right).
\end{split}
\end{equation*}
\end{proof}

By\eqref{212} and Sobolev's inequality, $\t\in L^2([0,T]; L^6(\O))\cap L^2([0,T]; L^3(\O))$,
since $a\ge 2$. Also, in view of the total energy conservation
inequality \eqref{21} and our assumptions, we know that $\r
P_e(\r)\in L^\infty([0,T]; L^1(\O))$ and $\r\t\in L^\infty([0,T];
L^1(\O))$. Thus, it is possible, by taking $\varepsilon$ small
enough, to get some \textit{ a priori } estimates from \eqref{22}
and \eqref{23} via Gronwall's inequality.

\begin{Lemma}[\bf{The control of $\int_\O\nabla
p\cdot\nabla\va(\r)\,dx$}]\label{l4}
\begin{equation*}
\begin{split}
-\int_\O\nabla p\cdot\nabla\va(\r)\,dx
&\le -\int_\O
\va'(\r)\t|\nabla\r|^2\,dx
 +\int_\O
\left(c_\varepsilon\kappa(\r,
\t)\f{|\nabla\t|^2}{\t^2}+\varepsilon\f{|\nabla\mu(\r)|^2}{\r}\right)\,dx \\
&\quad -\int_\O |\nabla\r^{-\f{l+1-\beta}{2}}|^2
\chi_{\{\r<A_1\}}\,dx,
\end{split}
\end{equation*}
with $A_1=\min\{A, A_0\}$.
\end{Lemma}
\begin{proof}
By our assumption \eqref{34}, we have
$$\nabla p=\t\nabla\r+\r\nabla\t+p_e'(\r)\nabla\r.$$
Hence, we have
\begin{equation}\label{214}
\begin{split}
-\int_\O\nabla p\cdot\nabla\va(\r)\,dx&=-\int_\O
\va'(\r)\t|\nabla\r|^2\,dx-\int_\O
\va'(\r)\r\nabla\t\cdot\nabla\r\,dx\\&\quad-\int_\O
p_e'(\r)\va'(\r)|\nabla\r|^2\,dx\\&\le-\int_\O
\va'(\r)\t|\nabla\r|^2\,dx-\int_\O
\va'(\r)\r\nabla\t\cdot\nabla\r\,dx\\&\quad-c\int_\O
|\nabla\r^{-\f{l+1-\beta}{2}}|^2 \chi_{\{\r<A_1\}}\,dx.
\end{split}
\end{equation}
because $\va'(\r)>0$ and $p_e'(\r)\ge 0$.

As for the second term on the right-hand side of \eqref{214}, we
have, by our assumption \eqref{33},
\begin{equation*}
\begin{split}
\left|\int_\O \va'(\r)\r\nabla\t\cdot\nabla\r\,dx\right|&\le
\int_\O |\va'(\r)\r\nabla\t\cdot\nabla\r|\,dx\\&\le \int_\O
\left(c_\varepsilon\kappa(\r,
\t)\f{|\nabla\t|^2}{\t^2}+\varepsilon\f{\r\t^2}{\kappa(\r,
\t)}\f{|\nabla\mu(\r)|^2}{\r}\right)\,dx\\&\le \int_\O
\left(c_\varepsilon\kappa(\r,
\t)\f{|\nabla\t|^2}{\t^2}+\varepsilon\f{|\nabla\mu(\r)|^2}{\r}\right)\,dx.
\end{split}
\end{equation*}
Thus, we have
\begin{equation*}
\begin{split}
-\int_\O\nabla p\cdot\nabla\va(\r)\,dx\le &-\int_\O
\va'(\r)\t|\nabla\r|^2\,dx+\int_\O \left(c_\varepsilon\kappa(\r,
\t)\f{|\nabla\t|^2}{\t^2}+\varepsilon\f{|\nabla\mu(\r)|^2}{\r}\right)\,dx\\&-\int_\O
|\nabla\r^{-\f{l+1-\beta}{2}}|^2 \chi_{\{\r<A_1\}}\,dx.
\end{split}
\end{equation*}
\end{proof}

Noting that Proposition \ref{p31} implies $\kappa(\r,
\t)\f{|\nabla\t|^2}{\t^2}\in L^1(\O\times(0,T))$. Therefore, it is
also possible, by incorporating the estimate into \eqref{22} and
\eqref{23}, to get some \textit{ a priori } estimates via
Gronwall's inequality.

\begin{Lemma}[\bf{The control of $\int_\O(\nabla\times\H)\times\H\cdot\f{\nabla\mu(\r)}{\r}\,dx$}]\label{l5}
$$\left|\int_\O(\nabla\times\H)\times\H\cdot\f{\nabla\mu(\r)}{\r}\,dx\right|\le c\int_\O\left(\f{|\nabla\times\H|^2\nu(\r, \t,
\H)}{\t}+\f{|\nabla\mu(\r)|^2}{\r}\right)\,dx.$$
\end{Lemma}
\begin{proof}
Indeed, we can estimate, by our assumption \eqref{36} and the
uniform bound of  $\H_n$ in $L^{\infty}(\O\times(0,T))$,
\begin{equation*}
\begin{split}
\left|\int_\O(\nabla\times\H)\times\H\cdot\f{\nabla\mu(\r)}{\r}\,dx\right|&
\le\int_\O|(\nabla\times\H)\times\H\cdot\f{\nabla\mu(\r)}{\r}|\,dx\\
&\le c\int_\O \left(\f{|\nabla\times\H|^2|\H|^2}{\r}+\f{|\nabla\mu(\r)|^2}{\r}\right)\,dx\\
&\le c\int_\O \left(\f{|\nabla\times\H|^2\nu(\r,
\t)}{\t}+\f{|\nabla\mu(\r)|^2}{\r}\right)\,dx.
\end{split}
\end{equation*}
\end{proof}

The entropy inequality \eqref{210} implies that
$\f{|\nabla\times\H|^2\nu(\r, \t)}{\t}$ belongs to
$L^1(\O\times(0,T))$, provided $\r_0 s_0$ and $\r_0\ln\r_0$ belong
to $L^1(\O)$. Therefore, from Lemma \ref{l1}-\ref{l5}, we can
deduce the following \textit{ a priori } estimates via Gronwall's inequality:
\begin{equation}\label{215}
\begin{cases}
&\|\sqrt{\r}\u\|_{L^\infty([0,T];L^2(\O))}\le c,\quad
\|\r^{-1/2}\nabla\mu(\r)\|_{L^\infty([0,T];L^2(\O))}\le c,\\
&\|(\r^{m/2}+\r^{\beta/2})\nabla\u\|_{L^2(\O\times(0,T))}\le
c,\quad
\|\f{\r^{m/2}+\r^{\beta/2}}{\t^{1/2}}\nabla\u\|_{L^2(\O\times(0,T))}\le
c,\\
&\|(\t)^{1/2}(\r^{(\beta-1)/2}+\r^{(m-1)/2})\f{\nabla\r}{\r}\|_{L^2(\O\times(0,T))}\le
c,\quad \|\r P_e(\r)\|_{L^\infty([0,T];L^1(\O))}\le c,\\
&\|\sqrt{\f{\r\t}{\mu'(\r)}}\nabla\va(\r)\|_{L^2(\O\times(0,T))}\le
c,\quad \|\r\t\|_{L^\infty([0,T];L^1(\O))}\le c,\\
&\|\H\|_{L^\infty([0,T];L^2(\O))}\le c,\quad
\|\sqrt{\nu}\nabla\times\H\|_{L^2(\O\times(0,T))}\le c,\\
&\|\sqrt{1+\r}\nabla\t^\alpha\|_{L^2(\O\times(0,T))}\le c, \quad
\|\nabla\r^{-\f{l+1-\beta}{2}}
\chi_{\{\r<A_1\}}\|_{L^2(\O\times(0,T))}\le c,
\end{cases}
\end{equation}
for all $\alpha\in[0, a/2]$.

\section{Some Integrability Lemmas}
As mentioned in \cite{hw2,l1}, the lack of \textit{ a priori }
estimates on approximation solutions to the compressible flow is
the main difficulty to prove the existence and the compactness of
global-in-time weak solutions. Indeed, the basic and natural \textit{ a
priori } estimates are not sufficient, since the energy equation
does not hold so far even in the distribution theory framework.
For more details, we refer the readers to \cite{hw2}.

This difficulty has been circumvented in \cite{f1, hw2} by
restricting the generality of the equations of state \eqref{33}
and \eqref{34}, and defining variational solutions for which the
energy equation \eqref{13} becomes two inequalities in the sense
of distributions. However, this approach requires significant
restrictions on the equations of state, in particular the ideal gas
case is excluded.

This section is devoted to the local integrability analysis of the
various energy fluxes such as $\r\u|\u|^2$, $\r\u e$, $\u p$,
$\kappa\nabla\t$. One of the crucial steps is the additional
integrability obtained on $\r$.

\subsection{Integrability of the velocity}

Let us begin with some bounds on the velocity with density dependent
weights.
\begin{Lemma}\label{l41}
Let $\O$ be either the whole space $\R^3$ or the three-dimensional
periodic box, and $T>0$. Let $\u$ be a vector field over
$\O\times(0,T)$ such that $\u\in L^{q_1}([0,T];
L_{loc}^{q_2}(\O))$, $\sqrt{\r}\u\in L^\infty([0,T]; L^2(\O))$,
and $\r\in L^\infty([0,T]; L_{loc}^p(\O))$ such that
$$q_1\in(1,2), \quad \textrm{and}\quad
\f{1}{p}+\f{2q_1}{q_2(q_1-1)}<1.$$ Then, there exists $\delta>3$
such that $\r^{1/3}\u\in L^\delta(\O'\times(0,T))$ for all bounded
subsets $\O'$ in $\O$.
\end{Lemma}
\begin{proof}
For the proof we refer the reader to Lemma 6.1 in \cite{bd1}.
\end{proof}

In order to apply Lemma \ref{l41} to improve the integrability of
$\r^{1/3}\u$, we need first to deduce the integrability of the
velocity $\u$. Indeed, following the computation in \cite{bd1},
one may write $\nabla\u=\r^{-\beta/2}\r^{\beta/2}\nabla\u$, and
then deduce that
\begin{equation}\label{a2}
\begin{split}
&\|\nabla\u\|_{L^{q_1}([0,T]; L^{q_3}(\O'))}\\
&\le C(\O')\left(1+\|\r^{-\beta/2}\chi_{\{\r<A_1\}}\|_{L^{2j}([0,T];
L^{6j}(\O'))}\right)
\|\r^{\beta/2}\nabla\u\|_{L^2(\O\times(0,T))}\\
&\le
C(\O')\left(1+\|\nabla\r^{-(l+1-\beta)/2}\chi_{\{\r<A_1\}}\|_{L^{2}([0,T];L^{2}(\O'))}\right)
\|\r^{\beta/2}\nabla\u\|_{L^2(\O\times(0,T))},
\end{split}
\end{equation}
with $$j=\f{l+1-\beta}{\beta},\quad
q_1=2\left(1-\f{\beta}{l+1}\right), \quad \textrm{and }
\f{1}{q_3}=\f{1}{6j}+\f{1}{2}.$$ In Lemma \ref{l41}, letting
$q_2=\f{3q_3}{3-q_3}=3q_1$, taking $p=6m-3$, we deduce that
$\r^{1/3}\u\in L^\delta((0,T)\times \O')$ for some $\delta>3$.

As a byproduct of previous analysis, we also can derive useful
bounds on energy flux. More precisely, one has for all bounded
subset $\O'$ of $\O$,
$$\|\r^{-l}\u\|_{L^s([0,T];L^r(\O'))}\le C\|\r^{-l}\|^{2/3}_{L^\infty([0,T];L^1(\O'))}\|\r^{-l}\|^{1/3}_{L^{2j_1}([0,T];L^{6j_1}(\O'))}
\|\u\|_{L^{q_1}([0,T];L^{q_2}(\O'))},$$ with
$$j_1=\f{l+1-\beta}{2l},\quad \f{1}{s}=\f{5l+3}{6(l+1-\beta)},\quad
\f{1}{r}=\f{17l+15-12\beta}{18(l+1-\beta)}.$$ Noting that the
hypothesis \eqref{35} implies that $s>1$ and $r>1$.

In order to bound the energy fluxes, it remains to control $\r^k\u$
in $L^\delta((0,T)\times \O')$ for some $\delta>1$. Since $\u$ is
bounded in $L^{q_1}([0,T]; L_{loc}^{q_2}(\O))$ and $\r$ is bounded
in $L^\infty([0,T]; L_{loc}^{(6m-3)/k}(\O))$, the hypothesis
\eqref{35} implies that $\r^k\u$ is bounded in $L^\delta((0,T)\times
\O')$ for some $\delta>1$.

\subsection{Integrability of the heat flux}
In this subsection, we will need the following integrability on
the temperature:
\begin{Lemma}\label{l42}
Let $\O$ be either the whole space $\R^3$ or a three-dimensional
periodic box and let $T>0$. Let $\t$ be a function over
$\O\times(0,T)$ such that $(\sqrt{\r}+1)\nabla\t^{a/2}$ and
$(\sqrt{\r}+1)\nabla\ln\t$ belong to $L^2(\O\times(0,T))$ with $a\ge
2$, $\r e\in L^\infty([0,T]; L^1(\O))$, and
$\r^{-1/2}\nabla\mu(\r)\in L^2(\O)$. Then, $\t^{\f{a-c+1}{2}}$
belongs to $L^2([0,T]; L^6(\O))$ for all $0<c\le 1$.
\end{Lemma}
\begin{proof}
For the proof, we refer the reader to Lemma 7.3 in \cite{bd1}.
\end{proof}

At this stage, in order to improve the integrability of the heat
flux, we need to use the following thermal energy equation (cf.
equation (1.15) in \cite{hw2}).
\begin{equation}\label{41}
\partial_t (\r \t)+\Dv(\r \t\u)-\Dv(\kappa(\r, \t)\nabla\t)=\nu|\nabla\times\H|^2+\Psi:\nabla\u-\t\r\Dv\u.
\end{equation}
As a matter of fact, we have
\begin{Lemma}\label{l43}
For any nondecreasing concave function from $R^+$ to $R$, one has
\begin{equation*}
\begin{split}
\int_\O f'(\t)&(2\mu
D(\u):D(\u)+\lambda|\Dv\u|^2+\nu|\nabla\times\H|^2)\,dx-\int_\O\kappa(\r,\t)f''(\t)|\nabla\t|^2\,dx\\
&\le\f{d}{dt}\int_\O\r f(\t)\,dx+\int_\O\r\t f'(\t)|\Dv\u|\,dx.
\end{split}
\end{equation*}
\end{Lemma}
\begin{proof}
Multiplying \eqref{41} by $f'(\t)$, one has
\begin{equation*}
\begin{split}
\int_\O f'(\t)&(2\mu
D(\u):D(\u)+\lambda|\Dv\u|^2+\nu|\nabla\times\H|^2)\,dx-\int_\O\kappa(\r,\t)f''(\t)|\nabla\t|^2\,dx\\
&=\int_\O f'(\t)(\t\r\Dv\u+\partial_t(\r\t)+\Dv(\r\t\u))\,dx\\
&=\int_\O f'(\t)(\t\r\Dv\u+\r\partial_t\t+\r\u\nabla\t)\,dx\\
&=\int_\O f'(\t)\t\r\Dv\u+\r\partial_tf(\t)+\r\u\nabla f(\t)\,dx\\
&=\int_\O f'(\t)\t\r\Dv\u+\partial_t(\r f(\t))+\Dv(\r\u f(\t))\,dx\\
&\le\int_\O f'(\t)\t\r|\Dv\u|\,dx+\f{d}{dt}\int_\O\r f(\t)\,dx,
\end{split}
\end{equation*}
here, we used twice the mass conservation equation \eqref{11}.
\end{proof}

Now, we consider $f'(\t)=\f{1}{\t^c}$ for some $0<c<1$ in Lemma
\ref{l43}, then we have
\begin{equation*}
\begin{split}
\int_\O \f{1}{\t^c}&(2\mu
D(\u):D(\u)+\lambda|\Dv\u|^2+\nu|\nabla\times\H|^2)\,dx+c\int_\O\kappa(\r,\t)\f{1}{\t^{c+1}}|\nabla\t|^2\,dx\\
&\le\f{d}{dt}\int_\O\r
\f{\t^{1-c}}{(1-c)}\,dx+\int_\O\r\t^{1-c}|\Dv\u|\,dx.
\end{split}
\end{equation*}

Keeping the hypothesis \eqref{33} in mind, we have
\begin{equation}\label{42}
\int_\O(1+\r)|\nabla(1+\t)^{(a-c+1)/2}|^2\,dx\le\f{d}{dt}\int_\O\r
\f{\t^{1-c}}{(1-c)}\,dx+\int_\O\r\t^{1-c}|\Dv\u|\,dx.
\end{equation}

For the second term on the right-hand side of \eqref{42}, one has
\begin{equation*}
\begin{split}
&\left|\int_\O
\r\t^{1-c}\Dv\u\,dx\right|\\
&\le\|\sqrt{3\lambda+2\mu}\Dv\u\|_{L^2}\\
&\quad\times(\|\r\t^{1-c}\chi_{\{\r<A\}}/\sqrt{3\lambda+2\mu}\|_{L^2}
+\|\r\t^{1-c}\chi_{\{\r\ge A\}}/\sqrt{3\lambda+2\mu}\|_{L^2})\\&\le
c\|\sqrt{3\lambda+2\mu}\Dv\u\|_{L^2}\\&\quad\times(\|\r^{2/5}\t^{1-c}\|_{L^2}\|
\r^{(6-5\beta)/10}\chi_{\{\r<A\}}\|_{L^\infty}
+A^{-m/2}\|\r\chi_{\{\r\ge A\}}\|_{L^6}\|\t^{1-c}\|_{L^3})\\&\le
c\|\sqrt{3\lambda+2\mu}\Dv\u\|_{L^2}\\
&\quad\times(\|\r\t^{1-c}\|^{2/5}_{L^1}\|
\t^{1-c}\|^{3/5}_{L^6}A^{(6-5\beta)/10}
+A^{(3-3m)/2}\|\r\chi_{\{\r\ge
A\}}\|^{m-1/2}_{L^{6m-3}}\|\t^{1-c}\|_{L^3}).
\end{split}
\end{equation*}
Thus, in view of Lemma \ref{l2} and Young's inequality, we have
\begin{equation*}
\begin{split}
&\left|\int_\O \r\t^{1-c}\Dv\u\,dx\right| \\
&\le
c\|\sqrt{3\lambda+2\mu}\Dv\u\|_{L^2}
\left(\|\r\t^{1-c}\|_{L^1}+\|\t^{1-c}\|_{L^6}+\|\t^{1-c}\|_{L^3}
\left\|\f{\nabla\mu(\r)}{\sqrt{\r}}\right\|_{L^2}\right)\\
&\le\varepsilon\|\sqrt{3\lambda+2\mu}\Dv\u\|_{L^2}^2
+c_\varepsilon\left(\|\r\t^{1-c}\|^2_{L^1}+\|\t^{1-c}\|^2_{L^6}+
\|\t^{1-c}\|^2_{L^3}
\left\|\f{\nabla\mu(\r)}{\sqrt{\r}}\right\|^2_{L^2(\O)}\right).
\end{split}
\end{equation*}
Because $\r\t$ and $\r$ belong to $L^\infty([0,T]; L^1(\O))$, we
deduce that $\r\t^{1-c}=\r^c(\r\t)^{1-c}$ belongs to
$L^\infty([0,T]; L^1(\O))$. Hence, integrating over t in the both
sides of \eqref{42}, and combining the estimates in Proposition
\ref{p31}, one has, if $0<c\le 1$,
\begin{equation}\label{43}
\int_0^T\!\!\!\!\int_\O(1+\r)|\nabla(1+\t)^{(a-c+1)/2}|^2\,dxdt\le
C,
\end{equation}
for some positive constant $C$. In particular, by Sobolev's
inequality, $\t^{a-c+1}$ belongs to $L^2([0,T]; L^3(\O))$.

\begin{Lemma}\label{l44}
Let $\O$ be either the whole space $\R^3$ or the three-dimensional
periodic box and let $T>0$. Assume that $\t$ satisfies the
conditions in Lemma \ref{l42} and
$\sqrt{1+\r}\nabla\t^{(a-c+1)/2}\in L^2([0,T];L^2(\O))$,
$\t^{a-c+1}\in L^2([0,T]; L^3(\O))$, $\r\in L^\infty([0,T];
L^{6m-3}(\O))$. Then, one has, for some $p>1$,
\begin{equation}\label{44}
\kappa(\r, \t)\nabla\t\in L^p([0,T]; L^p(\O)).
\end{equation}
\end{Lemma}
\begin{proof}
For the proof, we refer the reader to Lemma 6.2 in \cite{bd1}.
\end{proof}

\bigskip

\section{Compactness of Weak Solutions}
With the \textit{ a priori } estimates and integrability lemmas
obtained in the previous sections, we  now  study the compactness of
sequences of weak solutions $\{(\r_n, \u_n, \t_n,
\H_n)\}_{n=1}^{\infty}$ and pass to the limit in nonlinear terms.

To begin with, we will state the Aubin-Lions \textit{compactness
lemma} (see\cite{l3}, Ch. IV, and \cite{s1} for more recent
references) which we will use later. A simple statement goes as
follows:
\begin{Lemma}[\bf{Aubin-Lions Lemma}]\label{al}
Let $T>0$, $p\in(1,\infty)$ and let $\{f_n\}_{n=1}^{\infty}$ be a
bounded sequence of functions in $L^p([0,T]; H)$ where H is a
Banach space. If $\{f_n\}_{n=1}^{\infty}$ is also bounded in
$L^p([0,T]; V)$, where V is compactly imbedded in H and
$\{\partial f_n/\partial t\}_{n=1}^{\infty}$ is bounded in
$L^p([0,T]; Y)$ uniformly where $H\subset Y$, then
$\{f_n\}_{n=1}^{\infty}$ is relatively compact in $L^p([0,T]; H)$.
\end{Lemma}

\subsection{Compactness of the density}
From the uniform estimates derived in Lemma \ref{l2}, we deduce
that the sequence $\r_n$ is uniformly bounded in
$L^\infty([0,T];L^{6m-3}(\O'))$ for all bounded subset $\O'$ of
$\O$. Up to a subsequence, one may assume that $\r_n$ converges
weakly to some $\r$ in $L^2([0,T]; L^2_{loc}(\O))$. In fact, we
have
\begin{Lemma}\label{l51}
\begin{equation*}
\begin{split}
\partial_t(\r_n^m) \textrm{ is bounded in } L^2([0,T]; L^{3/2}(\O')),\\
\nabla\r_n^m \textrm{ is bounded in } L^\infty([0,T]; L^{3/2}(\O'));
\end{split}
\end{equation*}
as a consequence, up to a subsequence, $\r_n$ converges almost
everywhere and strongly to some element $\r$ in $C([0,T];
L^{p}(\O'))$ for all $1\le p<6m-3$. Moreover, the continuity
equation \eqref{11} hold in the sense of distributions.
\end{Lemma}
\begin{proof}
Let us consider the renormalized mass equation satisfied by
$h(\r)=\r^m$,
$$\partial_t(h(\r_n))+\Dv(h(\r_n)\u_n)+(m-1)\r_n^m\Dv\u_n=0.$$
The uniform bounds of $\sqrt{\r_n}\u_n$ in $L^\infty([0,T];
L^2(\O))$ and of $\r_n^{-1/2}h(\r_n)$ in $L^\infty([0,T]; L^6(\O'))$
imply that $h(\r_n)\u_n$ is bounded in $L^\infty([0,T];
L^{3/2}(\O'))$. On the other hand, $\r^{m/2}\Dv\u_n$ is bounded in
$L^2(\O\times(0,T))$, and the sequence $\r^{m/2}$ is bounded in
$L^\infty([0,T]; L^{6(2m-1)/m}(\O'))$, hence, $\r_n^m\Dv\u_n$ is
bounded in $L^2([0,T]; L^{3/2}(\O'))$. Thus, $\partial_t(h(\r_n))$
is uniformly bounded in $L^2([0,T]; L^{3/2}(\O'))$.

The spatial regularity of $\r_n^m$ can be estimated as follows.
Since $\nabla(\r_n^m)=\f{\nabla\r_n^m}{\sqrt{\r_n}}\sqrt{\r_n}$,
thus, $\nabla(\r_n^m)$ is bounded in $L^\infty([0,T];
L^{3/2}(\O'))$. The estimates deduced above, and thanks to
Aubin-Lions Lemma, give the strong convergence of $\r_n^m$ in
$L^2([0,T]; L^{2}(\O'))$. We denote the limit of $\r_n^m$ by
$\ov{\r^m}$. Then, since the function $h(s)$ is strictly increasing,
we conclude that $\r_n$ converges strongly to
$\r:=(\ov{\r^m})^{1/m}$ in $L^2([0,T]; L^{2}(\O'))$, and hence, by
interpolation, in $L^q([0,T]; L^{p}(\O'))$ for all $1\le p<6m-3$
with $\f{1}{p}=\f{1}{q}+\f{q-1}{(6m-3)q}$. And using the continuity
equation \eqref{11} again, we know actually $\r_n$ converges
strongly to $\r:=(\ov{\r^m})^{1/m}$ in $C([0,T]; L^{p}(\O'))$ for
all $1\le p<6m-3$.

Finally, we already know that
$\{\u_n\}_{n=1}^{\infty}$ is uniformly bounded in
$L^{q_1}([0,T]; W^{1,q_3}(\O'))$. Thus, up to a subsequence,
$\u_n$ converges weakly to some element $\u$ in $L^{q_1}([0,T];
W^{1,q_3}(\O'))$. By Sobolev's compact imbedding theorem, we also
have that $\u_n$ is uniformly bounded in $L^{5/3}([0,T]; L^{5}(\O'))$
since $q_1>5/3$ and $q_3>15/8$. As a consequence, we pass to the
limit in the mass conservation equation:
\begin{equation}\label{51}
\partial_t\r+\Dv(\r\u)=0, \quad \textrm{ in }
\mathcal{D}'(\O\times(0,T)).
\end{equation}
\end{proof}

In the spirit of Lemma \ref{l51}, we can pass to the limit in the
sense of distributions for the terms $\r_n P_e(\r_n)$, $p_e(\r_n)$,
$\r_n P_e(\r_n)\u_n$, $p_e(\r_n)\u_n$ since $\u_n$ converges weakly
to some element $\u$ in $L^{q_1}([0,T]; W^{1,q_3}(\O'))$.

\bigskip

\subsection{Compactness of the momentum}
In this subsection, we show the compactness of the momentum
$\r_n\u_n$ by Aubin-Lions Lemma. We already know from the previous
subsection that $\r_n\u_n$ converges weakly to $\r\u$ in
$L^2([0,T]; L^{3/2}(\O'))$, due to the facts $\sqrt{\r_n}\u_n\in
L^{\infty}([0,T]; L^2(\O))$ and $\r_n\in L^{\infty}([0,T];
L^3(\O'))$. To this end, we need to establish the uniform bound of
$\partial_t(\r_n\u_n)$ in some suitable functional space. Indeed,
we can show that the sequence $\partial_t(\r_n\u_n)$ is uniformly
bounded in $L^p([0,T];H^{-s}(\O))$ for some $p>1$ and $s$ large
enough.

From the momentum conservation equation \eqref{12}, we have
$$\partial_t(\r_n\u_n)=-\Dv\left(\r_n\u_n\otimes\u_n\right)-\nabla
p_n+(\na \times \H_n)\times \H_n+\Dv\Psi_n.$$
For the first term on the right-hand side, $\r_n\u_n\otimes\u_n$ is
bounded  in $L^{3/2}([0,T];L^{9/7}(\O'))$ uniformly as a product of
$\r_n^{1/3}\u_n$, bounded in $L^3([0,T]; L^3(\O'))$ and $\r_n$,
bounded in $L^\infty([0,T]; L^{6m-3}(\O'))$. For the second term,
$p_n=p(\r_n, \t_n)$ is bounded uniformly in $L^\infty([0,T];
L^{1}(\O))$, since $\r_n\t_n$ is uniformly bounded in
$L^{\infty}([0,T]; L^{1}(\O))$ and $p_e(\r_n)$ is bounded in
$L^\infty([0,T]; L^{1}(\O))$. For the third term, by \eqref{215},
$(\na \times \H_n)\times \H_n$ is bounded in
$L^2([0,T];L^{3/2}(\O))$. As for the fourth term $\Dv\Psi_n$, from
the fact $\sqrt{\mu(\r_n)}D(\u_n)$ and
$\sqrt{|\lambda(\r_n)|}\Dv\u_n$ are bounded in $L^2(\O\times(0,T))$,
and that $\sqrt{\mu(\r_n)}$ and $\sqrt{|\lambda(\r_n)|}$ are
uniformly bounded in $L^\infty([0,T]; L^6(\O'))$, we deduce that
$\Dv\Psi_n$ is bounded in $L^2([0,T]; L^{3/2}(\O'))$. Therefore, the
sequence $\partial_t(\r_n\u_n)$ is uniformly bounded in
$L^{3/2}([0,T];W^{-1,1}(\O'))$.

Next, we have
\begin{equation}\label{52}
\begin{split}
\partial_i(\r_n u_{nj})&=\r_n\partial_i
u_{nj}+u_{nj}\partial_i\r_n\\
&=\f{\r_n}{\r_n^{m/2}+\r_n^{\beta/2}}(\r_n^{m/2}+\r_n^{\beta/2})\partial_i
u_{nj}+\f{1}{\mu'(\r_n)}\sqrt{\r_n}u_{nj}\r_n^{-1/2}\partial_i\mu(\r_n).
\end{split}
\end{equation}
Using the hypothesis \eqref{32} and the estimate $\r_n\in
L^{\infty}([0,T]; L^{6m-3}(\O'))$, one can deduce that
$\f{\r_n}{\r_n^{m/2}+\r_n^{\beta/2}}\in L^{\infty}([0,T];
L^3(\O'))$ and $\f{1}{\mu'(\r_n)}\in L^{\infty}(\O\times(0,T))$.
Hence, from the estimates in \eqref{215}, we deduce that
$\r_n\u_n\in L^2([0,T]; W^{1,1}(\O'))$. Thus, by Aubin-Lions
Lemma, we deduce that $\r_n\u_n$ converges strongly to $\r\u$ in
$L^{3/2}([0,T]; L^p(\O'))$ for all $1\le p<3/2$.

As a conclusion, the product $\r_n|\u_n|^2$, converges strongly to
$\r|\u|^2$ in $L^1(\O\times(0,T))$, since $\r_n\u_n$ converges
weakly to $\r\u$ in $L^{\infty}([0,T]; L^{3/2}(\O'))$, strongly to
$\r\u$ in $L^{3/2}([0,T]; L^p(\O'))$ for all $1\le p<3/2$ and $\u_n$
is bounded uniformly in $L^{5/3}([0,T]; L^5(\O'))$. Using the fact
$\r_n^{1/3}\u_n=\r_n^{1/3}\u_n \chi_{\{\r_n\le
\varepsilon\}}+\r_n^{1/2}\u_n \r_n^{-1/6}\chi_{\{\r_n>
\varepsilon\}}$, $\r_n^{1/3}\u_n$ is the sum of a uniformly small
term in $L^1(\O'\times(0,T))$ and another term converging to
$\r^{1/3}\u \chi_{\{\r>\varepsilon\}}$ and then we deduce that
$\r_n^{1/3}\u_n$ converges strongly in $L^1([0,T]; L^1(\O'))$ to
$\r^{1/3}\u$. Finally, by the interpolation and the uniform bound of
$\r_n^{1/3}\u_n$ in $L^\delta(\O'\times(0,T))$ for some $\delta>3$,
we conclude that $\r_n^{1/3}\u_n$ converges strongly in
$L^3(\O'\times(0,T))$ to $\r^{1/3}\u$.

\bigskip

\subsection{Compactness of the temperature}
In this subsection, we want to derive compactness results for the
 energy $\E_n$ and the temperature $\t_n$. The first step
is to derive uniform bounds in $L^p([0,T];H^{-s}(\O))$ for some
$p>1$ and $s$ large enough for the sequence $\partial_t(\E_n)$.
Indeed, we can rewrite the energy conservation equation
\eqref{13}:
$$\partial_t(\E_n)=-\Dv(\u_n(\E_n'+p_n))+\Dv((\u_n\times\H_n)\times\H_n+\nu_n\H_n\times(\nabla\times\H_n)+\u_n\Psi_n+\kappa_n\nabla\theta_n).$$

For the first term of the right-hand side, we
already know that $\r_n\u_n|\u_n|^2$ is uniformly bounded in
$L^q([0,T]; L^q(\O'))$ for some $q>1$. Also, we already get the
uniform bounds in $L^q([0,T]; L^q(\O'))$ with some $q>1$ for the
terms $\r_n^{-l}\u_n$, $\r_n^k\u_n$ and
$\kappa(\r_n,\t_n)\nabla\t_n$ in  Section 4. Next, the uniform
bound of $\r_n\u_n$ in $L^{\infty}([0,T]; L^{3/2}(\O'))$ and the
uniform bound of $\t_n$ in $L^2([0,T]; L^6(\O))$ implies that
$\r_n\u_n\t_n$ is bounded in $L^q([0,T]; L^q(\O'))$ for some
$q>1$. Hence, $\r_n\u_n e_n$ and $p_n\u_n$ are bounded in
$L^q([0,T]; L^q(\O'))$ for some $q>1$.

For the viscous flux $\u_n\Psi_n$, we note the facts that
$$\sqrt{\mu(\r_n)}D(\u_n) \text{ and }\sqrt{|\lambda(\r_n)|}\Dv\u_n \text{ are
bounded in } L^2(\O\times(0,T)),$$ and
$$\sqrt{\mu(\r_n)}\r_n^{-1/3} \text{ and }
\sqrt{|\lambda(\r_n)|}\r_n^{-1/3} \text{ are bounded in }
L^\infty([0,T];  L^{18(2m-1)/(3m-2)}(\O')) $$ due to the hypothesis
\eqref{32} and Lemma \ref{l2}, hence in $L^\infty([0,T];
L^{9}(\O'))$, and $\r_n^{1/3}\u_n$ is bounded in
$L^3(\O\times(0,T))$. Thus,  the viscous fluxes are
bounded in $L^{6/5}([0,T]; L^{18/17}(\O'))$.

As for the terms related to the magnetic field, we have the fact
$\H_n$ is bounded in $L^\infty([0,T]; L^2(\O))\cap L^2([0,T],
H^1(\O))$, hence, by interpolation, in $L^p([0,T];L^q(\O))$ where
$p>5$, $q<30/11$, and $$\f{1}{q}=\f{1}{2}-\f{2}{3p}.$$ Thus,
$(\u_n\times\H_n)\times\H_n$ belongs to $L^q([0,T]; L^q(\O'))$ for
some $q>1$, since $\u_n$ is uniformly bounded in $L^{5/3}([0,T];
L^{5}(\O'))$. By the bound of the magnetic field coefficient
$\nu(\r, \t)$, we know $\nu_n\H_n\times(\nabla\times\H_n)$ belongs
to $L^{2p/(2+p)}([0,T];L^{2q/(2+q)}(\O'))$ for $p>5$, $2<q<30/11$,
hence, in $L^q([0,T]; L^q(\O'))$ for some $q>1$.

With the bound of $\partial_t(\E_n)$ in mind, we can show the
strong convergence of the term $\r_n\t_n^2$. To this end, we will
follow the argument in \cite{bd1}. First, we note that the strong
convergence of $\sqrt{\r_n}\u_n$ to $\sqrt{\r}\u$ in $L^r([0,T];
L^2(\O'))$ and the strong convergence of $\r_n P_e(\r_n)$ to $\r
P_e(\r)$ in $L^r([0,T]; L^1(\O'))$ for all $r\in (1,\infty)$. Let
us introduce
$$\K=\{f\in L^|_{loc}(\O)|\|\nabla f\|_{L^2(\O)}=1\}$$
and a sequence $\T_k$ of regularizing kernels given for instance by
convolution operators such that the following basic properties hold:
$$\sup_{f\in\K}\|f-\T_k f\|\le \f{C}{k},$$
and for all compact subset $\O'\subset\O$, there exists $C_{k,\O'}$
such that for all $f\in \K$,
$$\|\T_k f\|_{L^\infty(\O')}\le C_{k,\O'} \quad \textrm{ and } \; \T_k
f\in H^s(\O) \quad \textrm{ for all } s>0.$$ Then, we can deduce
that for any compact subset $\O'\subset\O$, one has
\begin{equation}\label{52}
\begin{split}
\left|\int_{\O'\times(0,T)}(\r_n\t_n^2-\r\t^2)\,dx\,dt\right|&\le
\f{C}{k}(\|\r_n\t_n\|_{L^2(\O'\times(0,T))}+\|\r\t\|_{L^2(\O'\times(0,T))})\\
&\quad+\left|\int_{\O'\times(0,T)}\r\t(\t_n-\t)\,dx\,dt\right|\\
&\quad+\f{1}{2}\|(|\H_n|^2-|\H|^2)\T_k\t_n\|_{L^1(\O'\times(0,T))}\\
&\quad+\f{1}{2}\|(\r_n|\u_n|^2-\r|\u|^2)\T_k\t_n\|_{L^1(\O'\times(0,T))}\\
&\quad+\|(\r_n P_e(\r_n)-\r
P_e(\r))\T_k\t_n\|_{L^1(\O'\times(0,T))}\\
&\quad+\left|\int_{\O'\times(0,T)}(\E_n-\E)\T_k\t_n\,dx\,t\right|.
\end{split}
\end{equation}
Let us observe that the first term of the above right-hand side is
bounded  by
$$\f{C}{k}(\|\r_n\|_{L^\infty([0,T];L^3(\O'))}\|\t_n\|_{L^2([0,T];L^6(\O'))}
+\|\r\|_{L^\infty([0,T];L^3(\O'))}\|\t\|_{L^2([0,T];L^6(\O'))}).$$
Therefore, given $\varepsilon>0$, there exists an integer $k_0$
such that the preceding term is less than $\varepsilon/4$
uniformly in $n$. Dealing with the second term is an easy task
since $\t_n$ converges weakly to $\t$ in $L^2([0,T];L^6(\O'))$, so
that for $n$ large enough, the second term is estimated by
$\varepsilon/4$. But $\|\T_{k_0}\t_n\|_{L^\infty(\O')}$ is
uniformly bounded in $n$, whereas $\H_n^2$, $\r_n|\u_n|^2$ and
$\r_n P_e(\r_n)$ converges strongly respectively to $\H^2$,
$\r|\u|^2$ and $\r P_e(\r)$ in $L^1(\O'\times(0,T))$ respectively,
so that the sum of the third and the fourth term is estimated by
$\varepsilon/4$ for large enough $n$. For the last term with
$k=k_0$, the uniform bound of $\partial_t(\E_n)$ in
$L^p([0,T];H^{-s}(\O))$ where $p>1$ implies that up to a
subsequence, $\E_n$ converges strongly to $\E$ in
$C([0,T];H^{-s}(\O))$, so that, for large enough $n$, the
right-hand side of above inequality is less than $\varepsilon$. It
follows that $\sqrt{\r_n}\t_n$ converges strongly in
$L^2(\O\times(0,T))$.

On the other hand, we know the strong convergence of $\r_n^{-1/2}$
to $\r^{-1/2}$ in $L^2([0,T];L^{p}(\O'))$ for $p<6$ due to the
Lebesgue's dominated convergence theorem, the estimate \eqref{215}
and the strong convergence of the density. And then, we deduce that
$\t_n$ converges to $\t$ in $L^1([0,T];L^r_{loc}(\O'))$ for all
$r<3/2$. Recalling the uniform bound of $\t_n^{a/2}$ in
$L^2([0,T];L^{6}(\O'))$, we deduce that $\t_n$ converges strongly to
$\t$ in $L^p([0,T];L^{q}(\O'))$ for all $p<a$, and $q<3a$ with
$$\f{1}{q}=\f{a-p}{p(a-1)r}+\f{p-1}{3p(a-1)},$$ for all $r<3/2$.

\subsection{Compactness of the magnetic field}

The aim of this subsection is to show the compactness of the
magnetic field. From Lemma \ref{l1}, and the hypothesis \eqref{36},
we deduce that $\H_n\in L^2([0,T]; H^1(\O))\cap L ^{\infty}([0,T];
L^2(\O))$. Thus, we can assume that $\H_n$ converges weakly to some
element $\H$ with $\Dv\H=0$ in $L^2([0,T]; H^1(\O))\cap L
^{\infty}([0,T]; L^2(\O))$.

On the other hand, from the equation \eqref{14}, we know that
\begin{equation}\label{531}
\partial_t\H_n=\nabla\times(\u_n\times\H_n)-\nabla\times(\nu(\r_n,
\t_n)\nabla\times\H_n). \end{equation} For the first term on the
right-hand side of \eqref{531}, we deduce that $$\u_n\times\H_n\in
L^{5/3}([0,T]; L^{10/7}(\O'))$$ because of $\u_n \in
L^{5/3}([0,T]; L^{5}(\O'))$ and $\H_n \in L^{\infty}([0,T];
L^{2}(\O))$. For the second term on the right-hand side of
\eqref{531}, we have $\nu(\r_n, \t_n)\nabla\times\H_n \in
L^2(\O\times(0,T))$ due to the hypothesis \eqref{36}. Hence,
$\partial_t\H_n$ is bounded in $L^{5/3}([0,T]; W^{-1,10/7}(\O'))$.
By Aubin-Lions Lemma, we deduce that $\H_n$ converges strongly to
$\H$ in $L^{5/3}([0,T]; L^p(\O'))$ for any $5<p<6$. Furthermore,
due to the uniform bound of $\H_n$ in $L^{\infty}([0,T];
L^2(\O))$, by using the interpolation, one obtain that $\H_n$
converges strongly to $\H$ in $L^{p}([0,T]; L^q(\O'))$ for some
$p>5$ and some $q>5/2$.

\subsection{Proof of Theorem \ref{mt}}
To finish the proof of Theorem \ref{mt}, we need to check that the
limit functions $\r$, $\u$, $\t$, $\H$ are indeed the weak
solutions, as defined in the introduction. We will complete this
proof by several steps.

\textit{Step 1: Convergence of the mass conservation equation.}\quad
 Let
us start with the mass conservation equation \eqref{11}, since
$\r_n$ converges strongly to $\r$ in $C([0,T]; L^p(\O'))$ for all
$1\le p<6m-3$, and $\u_n$ converges weakly to $\u$ in
$L^{q_1}([0,T]; W^{1,q_3}(\O))$, we deduce that, by the Sobolev's
compact imbedding theorem, $\r_n\u_n$ converges strongly to $\r\u$
in $L^r([0,T]; L^1(\O'))$ for some $r>1$. In particular, the mass
conservation equation \eqref{11} is satisfied in the sense of
distributions.

\textit{Step 2: Convergence of the momentum conservation equation.} \quad
For the momentum conservation equation \eqref{12}, the strong
convergence of $\r_n\u_n$ and $\r_n\u_n\otimes\u_n$ in $L^1([0,T];
L^1(\O'))$ can ensure the passing to limit in the sense of
distribution for the two corresponding term in the momentum
conservation equation \eqref{12}. On the other hand, since $\H_n$
converges weakly* to $\H$ in $L^\infty([0,T]; L^2(\O))\cap
L^2([0,T]; H^1(\O))$, this implies that the nonlinear term
$(\nabla\times\H_n)\times\H_n$ converges to
$(\nabla\times\H)\times\H$ in the sense of distributions. As a
product of $\r_n$ and $\t_n$, which respectively converge strongly
in $C([0,T]; L^2(\O'))$ and in $L^2(\O\times(0,T))$, the term
$\nabla(\r_n\t_n)$ converges to the limit $\nabla(\r\t)$ in the
sense of distributions. The term $p_e(\r)$ is already done in view
of the hypothesis \eqref{35} and the strong convergence of $\r$ in
$C([0,T]; L^p(\O'))$ for all $1\le p<6m-3$. Thus, we are left to
show the convergence of the viscous flux. In fact,
\begin{equation}\label{53}
\mu(\r_n)D(\u_n)=D(\mu(\r_n)\u_n)-\f{1}{2}\left(\sqrt{\r_n}\u_n
\otimes\f{\nabla\mu(\r_n)}{\sqrt{\r_n}}+
\f{\nabla\mu(\r_n)}{\sqrt{\r_n}}\otimes\sqrt{\r_n}\u_n\right).
\end{equation}
Since $\f{\mu(\r_n)}{\sqrt{\r_n}}$ converges strongly to
$\f{\mu(\r)}{\sqrt{\r}}$ in $L^\infty([0,T]; L^2(\O'))$ and
$\sqrt{\r_n}\u_n$ converges strongly to $\sqrt{\r}\u$ in $L^2([0,T];
L^2(\O'))$, the fist term on the right-hand side of \eqref{53}
converges to the corresponding term in the sense of distributions.
The convergence of the second term on the right-hand side of
\eqref{53} in the sense of distributions can be shown by using the
weak convergence of $\r_n^{-1/2}\nabla\mu(\r_n)$ to
$\r^{-1/2}\nabla\mu(\r)$ in $L^2([0,T]; L^2(\O'))$ and the strong
convergence $\sqrt{\r_n}\u_n$ in $L^2([0,T]; L^2(\O'))$. For the
bulk viscous term $\lambda(\r_n)\Dv\u_n$, by the assumption
\eqref{31}, it may be written in the renormalized sense:
$$\lambda(\r_n)\Dv\u_n=-2(\partial_t\mu(\r_n)+\Dv(\mu(\r_n)\u_n)),$$
which can be shown directly by the convergence of $\r_n$ and $\u_n$,
and hence, the convergence in the sense of distributions for the
momentum conservation equation is done.

\textit{Step 3: Convergence of the energy conservation equation.} \quad
The
main difficulties in this step lie in the passage to the limit for
the energy flux $\u(\E'+p)$, the heat flux $\kappa\nabla\t$, the
viscous term $\u\Psi$, and the nonlinear terms
$(\u\times\H)\times\H$, $\nu\H\times(\nabla\times\H)$, because we
already showed that $\E_n$ converges strongly to $\E$ in $C([0,T];
H^{-s}(\O))$ for some $s>0$.

For the energy flux $\r_n\u_n\t_n$, since $\sqrt{\r_n}\u_n$ and
$\sqrt{\r_n}\t_n$ converge strongly in $L^2(\O'\times(0,T))$ to
$\sqrt{\r}\u$ and $\sqrt{\r}\t$ respectively, $\r_n\u_n\t_n$
converges strongly in $L^1(\O'\times(0,T))$ to $\r\u\t$. For the
energy flux $\r_n\u_n|\u_n|^2$, the strong convergence of
$\r_n^{-1/2}$ in $C([0,T]; L^p(\O'))$ for all $p<6$ implies that
$\r_n^{-1/6}$ converges strongly to $\r^{-1/6}$ in $C([0,T];
L^3(\O'))$. Hence, the term $\r_n^{-1/6}\sqrt{\r_n}\u_n$ converges
strongly to $\r^{-1/6}\sqrt{\r}\u$ in $L^2([0,T]; L^{6/5}(\O'))$,
because of the strong convergence of $\sqrt{\r}\u$ in
$L^2(\O'\times(0,T))$. And Lemma \ref{l41} implies that
$\r_n^{1/3}\u_n$ is uniformly bounded in $L^\delta(\O'\times(0,T))$
for some $\delta>3$. This fact, combining with the interpolation
inequality, gives the strong convergence of $\r\u|\u|^2$ in
$L^1(\O'\times(0,T))$. The analysis at the end of Section 4 tells
the strong convergence of $\r_n\u_n P_e(\r_n)$ and $\u_n p_n$ in
$L^1(\O'\times(0,T))$ to $\r\u P_e(\r)$ and $\u p$, respectively.
Thus, the energy flux $\u_n(\E'_n+p_n)$ converges strongly to
$\u(\E'+p)$ in $L^1(\O'\times(0,T))$.

The strong convergence of $\t$ in $L^p([0,T];L^{q}(\O'))$ for $p<a$
and $q<3a$ implies that $\t_n^{a/2}$ converges strongly in
$L^{2}([0,T]; L^3(\O'))$ to $\t^{a/2}$. This fact, together the
strong convergence of $\r_n$ in $C([0,T]; L^p(\O'))$ for $p<6m-3$,
implies that $\sqrt{1+\r_n}(1+\t_n^{a/2})$ converges to
$\sqrt{1+\r}(1+\t^{a/2})$ in $L^{2}(\O'\times(0,T))$. That means
$\kappa^{1/2}(\r_n, \t_n)$ strongly converges to $\kappa^{1/2}(\r,
\t)$ in $L^{2}(\O'\times(0,T))$. Similarly, it follows that
$(1+\r_n)^{1/2}\t_n^{(a+c+1)/2}$ converges strongly to
$(1+\r)^{1/2}\t^{(a+c+1)/2}$ in $L^2(\O'\times(0,T))$ due to the
strong convergence of $\r_n$ and $\t_n$. Therefore, $\kappa_0(\r_n,
\t_n)(1+\r_n)^{1/2}(1+\t_n)^{(a+c+1)/2}$ converges strongly to
$\kappa_0(\r, \t)(1+\r)^{1/2}(1+\t)^{(a+c+1)/2}$ in
$L^2(\O'\times(0,T))$. On the other hand, we deduce from \eqref{43}
that $(1+\r_n)^{1/2}\nabla(1+\t_n)^{(a-c+1)/2}$ is uniformly bounded
in $L^2(\O\times(0,T))$, hence weakly converges to some element
$\omega$ in $L^2(\O\times(0,T))$. It also follows that
$\nabla(1+\t_n)^{(a-c+1)/2}$ is uniformly bounded in
$L^2(\O\times(0,T))$, and hence weakly converges to
$\nabla(1+\t)^{(a-c+1)/2}$ in $L^2(\O'\times(0,T))$. Due to the
strong convergence of $\r_n$ in $L^2(\O'\times(0,T))$, we deduce
that $\omega=(1+\r)^{1/2}\nabla(1+\t)^{(a-c+1)/2}$. Finally, we
write
\begin{equation*}
\begin{split}
\kappa(\r_n, \t_n)\nabla\t_n&=\kappa_0(\r_n,
\t_n)(1+\r_n)(1+\t_n)^a\nabla\t_n\\&=\kappa_0(\r_n,
\t_n)(1+\r_n)^{1/2}(1+\t_n)^{(a+c+1)/2}(1+\r_n)^{1/2}\nabla(1+\t_n)^{(a-c+1)/2}.
\end{split}
\end{equation*}
This, together with the strong convergence of $\kappa_0(\r_n,
\t_n)(1+\r_n)^{1/2}(1+\t_n)^{(a+c+1)/2}$ and the weak convergence
of $(1+\r_n)^{1/2}\nabla(1+\t_n)^{(a-c+1)/2}$, implies that
$\kappa(\r_n,\t_n)\nabla\t_n$ converges to $\kappa(\r,\t)\nabla\t$
at least in the sense of distributions.

For the viscous terms, $\sqrt{\mu(\r_n)}D(\u_n)$ and
$\sqrt{|\lambda(\r_n)|}\Dv\u_n$ converges weakly to
$\sqrt{\mu(\r)}D(\u)$ and $\sqrt{|\lambda(\r)|}\Dv\u$ respectively
in $L^2(\O'\times(0,T))$, because of the hypothesis \eqref{32},
the uniform bound on $\r_n$ in Lemma \ref{l2}, and the uniform
estimate \eqref{a2}. On the other hand, $\r_n^{1/3}\u_n$ strongly
converges to $\r^{1/3}\u$ in $L^3(\O'\times(0,T))$ and
$\r_n^{-1/3}\sqrt{\mu(\r_n)}$, as well as
$\r_n^{-1/3}\sqrt{\lambda(\r_n)}$ converges strongly to
$\r^{-1/3}\sqrt{\mu(\r)}$ and $\r^{-1/3}\sqrt{\lambda(\r)}$
respectively in $L^{\infty}([0,T]; L^6(\O'))$. Hence $\Psi_n\u_n$
converges to $\Psi\u$ at least in the sense of distributions.

Finally, we deal with the convergence of two nonlinear terms:
$(\u_n\times\H_n)\times\H_n$ and $\nu(\r_n,
\t_n)\H_n\times(\nabla\times\H_n)$. First, since $\H_n$ converges
strongly to $\H$ in $L^{p}([0,T]; L^q(\O'))$ for some $p>5$ and
some $q>5/2$, $\u_n\times\H_n$ weakly converges to $\u\times\H$ in
$L^{p}([0,T]; L^q(\O'))$ for some $p>5/4$ and $q>5/3$ because
$\u_n$ converges weakly to $\u$ in $L^{5/3}([0,T]; L^5(\O'))$.
From this, we can deduce that $(\u_n\times\H_n)\times\H_n$
converges to $(\u\times\H)\times\H$ in the sense of distributions.
Second, the strong convergence of $\r_n$, $\t_n$, $\H_n$ and the
hypothesis \eqref{36} imply that $\nu(\r_n,\t_n)\H_n$ converges
strongly to $\nu(\r,\t)\H$ in $L^{p}([0,T]; L^q(\O'))$ for some
$p>5$ and some $q>5/2$. By the weak convergence of $\H_n$ in
$L^2([0,T]; H^1(\O))$, one deduce that $\nu(\r_n,
\t_n)\H_n\times(\nabla\times\H_n)$ converges to $\nu(\r,
\t)\H\times(\nabla\times\H)$ at least in the sense of
distributions.

Therefore, the energy conservation equation \eqref{13} holds at
least in the sense of distributions.

\textit{Step 4: Convergence of the magnetic field equation.} \quad
Similar to the argument in Step 3, we can show that $\nu(\r_n,
\t_n)\nabla\times\H_n$ converges weakly to
$\nu(\r,\t)\nabla\times\H$ in $L^2(\O'\times(0,T))$. Also the strong
convergence of $\H_n$ in $L^{p}([0,T]; L^q(\O'))$ for some $p>5$ and
some $q>5/2$ and the weak convergence of $\u_n$ in $L^{5/3}([0,T];
L^5(\O'))$, imply that $\u_n\times\H_n$ converges to $\u\times\H$ at
least in the sense of distributions. Hence, the induction equation
holds at least in the sense of distributions.

The proof is complete.

\bigskip\bigskip

\section*{Acknowledgments}

Xianpeng Hu's research was supported in part by the National Science
Foundation grant DMS-0604362. Dehua Wang's research was supported in
part by the National Science Foundation grants DMS-0244487,
DMS-0604362, and the Office of Naval Research grant
N00014-07-1-0668.

\end{document}